\documentclass[12pt,a4paper]{article}

\usepackage{paralist}
\usepackage{graphics}
\usepackage{graphicx}
\usepackage{epsfig}
\usepackage[latin1]{inputenc}
\usepackage[T1]{fontenc}
\usepackage{amsmath}
\usepackage{amsfonts}
\usepackage{ae}
\usepackage{units}
\usepackage{icomma}
\usepackage{color}
\usepackage{graphicx}

\begin{document}

\begin{center}
\Large \textbf{Reconstruction of dielectric constants of multi-layered optical fibers
using propagation constants measurements}
\end{center}

\begin{center}
E. M. Karchevskii,$^1$ L. Beilina$^2$, A. O. Spiridonov,$^1$ and A. I. Repina$^1$
\end{center}

\noindent $^1$Kazan Federal University, Kremlevskaya 18, Kazan 420008, Russia

\noindent $^2$Chalmers University of Technology and University of Gothenburg, 42196
Gothenburg, Sweden

\begin{abstract}

We present new method for the numerical reconstruction of the variable
refractive index of multi-layered circular weakly guiding dielectric
waveguides using the measurements of the propagation constants of
their eigenwaves.  Our numerical examples show stable reconstruction
of the dielectric permittivity function $\varepsilon$ for random noise
level using these measurements.

\end{abstract}

\section{Introduction}

The development of analytical methods for the study of mathematical
models, the development and the justification of efficient numerical
methods for the solution of spectral problems of the theory of
dielectric waveguides attract the much attention (see, for example,
\cite{BKarch, Coldren_2012, Obayya_2011}). Dielectric waveguides form
the basic components of the microdevices used in the field of
integrated optics, photonics, and laser technology. The development of
analytical methods for the study of mathematical models of optical
microdevices, the development of numerical methods for accurate and
stable computations of their characteristics are essential for
designing and optimizing of such devices.

The problem of reconstruction of the refractive index of an
inhomogeneous dielectric waveguide from the measurements of the
propagation constants of its eigenwaves is very urgent.  Universal
numerical methods of reconstruction of refractive index of dielectric
objects are designed for coefficient inverse diffraction problems and
ignore the waveguide properties of the devices (see, for example,
\cite{BOOK}). The methods for the determination of the optical
characteristics of dielectric waveguides are proposed for waveguides
of some special forms (see, for example, \cite{Sokolov_2015,
  Sotsky_2015}).  For instance, the waveguide spectroscopy is widely
used for planar (one-dimensional) multi-layered waveguides
\cite{Khomchenko_2005}.  For such waveguides the characteristic
equation (a transcendental equation which connects the refractive
indices of the waveguide's layers with the propagation constants of
its eigenwaves) is well known. The method consists in minimization of
a functional, depending on the refractive indices of the layers. The
value of the functional is equal to the distance between the vector of
calculated (as the roots of the characteristic equation) propagation
constants and the vector of experimentally measured propagation
constants. Considerable efforts have been directed towards the
development of effective methods of minimization of this functional
\cite{Robert_2007, Schneider_2008}.

This method was extended in our previous works to the two-dimensional problem for the
waveguide with the piecewise-constant refractive index and an arbitrary cross-sectional
boundary \cite{Karchecvskii_Spiridonov_Beilina_SPMS_2015,
Karchevskii_Spiridonov_Repina_Beilina_PRI_2015}. The role of the ``characteristic
equation" is played by the system of two boundary weakly singular integral equations
whose kernels depend on the propagation constants and refractive indices of the vaweguide
and the environment. We presented new numerical methods for the solution of the inverse
spectral problem to determine the dielectric constants of core and cladding in optical
fibers. These methods use measurements of propagation constants. Our algorithms are based
on approximate solution of a nonlinear nonselfadjoint eigenvalue problem for a system of
weakly singular integral equations. We studied the inverse problem and proved that it is
well posed. Our numerical results indicated good accuracy of new algorithms.  Clearly,
the generalization of this method for an inhomogeneous waveguide is an urgent task.
Theoretical justification of this method is a very interesting problem.

In the presented paper we construct a Tikhonov functional and propose a method for the
numerical reconstruction of the variable refractive index of multi-layered circular
weakly guiding dielectric waveguides from the measurements of the propagation constants
of their eigenwaves. We theoretically investigate the Tikhonov functional and demonstrate
the practical effectiveness of the proposed method.  The success in the Tikhonov
regularization method largely depends on a good initial approximation to the solution.
The matematical analysis of the forward problem has allowed us to find a good first
approximation to the refractive index of the waveguide \cite{Frolov_Kartchevskiy_2013,
  Spiridonov_Karchevskiy_DD_2013}. This is confirmed by the numerical
experiments of solving the similar inverse problems
\cite{Karchecvskii_Spiridonov_Beilina_SPMS_2015,
  Karchevskii_Spiridonov_Repina_Beilina_PRI_2015}. Using this initial
approximation, in this paper we compute solutions of the investigated
problem and demonstrate that the algorithm is accurate and stable.

\section{Main equations}

\begin{figure}[h!]\centering
  \includegraphics[width=0.65\textwidth]{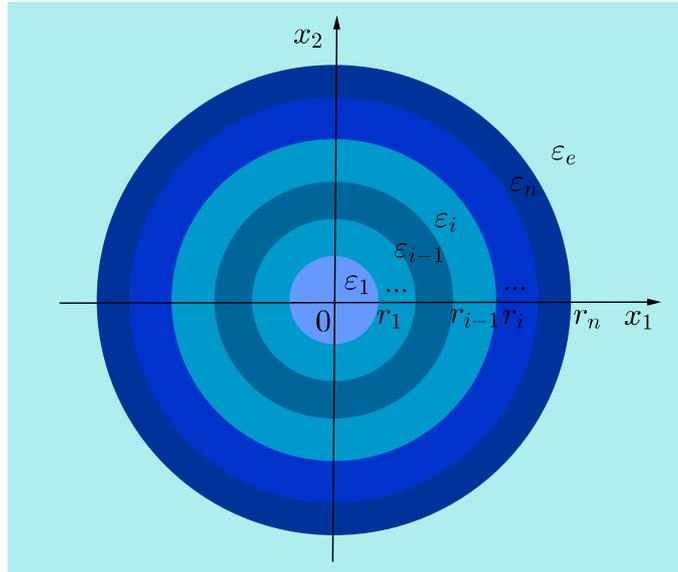}
  \caption{Geometry of a  multi-layered circular optical fiber}\label{plot_1}
\end{figure}

Let us consider a  multi-layered circular optical fiber shown at Figure ~\ref{plot_1} as
a regular cylindrical dielectric waveguide in a free space. The axis of the cylinder is
parallel to the $x_3$-axis. Suppose that the circles separating the layers of the
waveguide have the radiuses $r_l$, $l=1,\,2,\,\ldots,\, n$. Let the permittivity be
prescribed as a positive piecewise constant function~$\varepsilon$ which is equal to a
constant~$\varepsilon_e$ for $r=\sqrt{x_1^2+x_2^2}>r_n$ and to constants $\varepsilon_l$,
where $l=1,\,2,\,\ldots,\, n$, in corresponding layers. Let us denote
$\varepsilon_+=\max\limits_{l=1,\,2,\,\ldots,\, n}\varepsilon_l$, and suppose that
\begin{equation}\label{kem_perm_cond_1}
\min\limits_{l=1,\ldots,n}\varepsilon_l > \varepsilon_e \geq 1.
\end{equation}

Eigenvalue problems of optical waveguide theory~\cite{Snider_Love_1983} are formulated on
the base of the set of homogeneous Maxwell equations
\begin{equation} \label{Maxvell1}
\begin{array}{c}
\displaystyle {\textrm{rot} \cal E  = } - \mu_0 \frac{{\partial {\cal H}}}{{\partial t}},
\quad
 {\textrm{rot} \cal H  = }\varepsilon_0 \varepsilon
\frac{{\partial {\cal E}}}{{\partial t}}.
\end{array}
\end{equation}
Here, $\cal E$ and $\cal H$ are  electric and magnetic field vectors;  $\varepsilon_0$
and $\mu_0$ are the free-space dielectric and magnetic constants. Nontrivial solutions of
set (\ref{Maxvell1}) which have the form
\begin{equation} \label{modes}
\begin{array}{c}
\displaystyle \left[ {\begin{array}{*{20}c}
   {\cal E}  \\
   {\cal H}  \\
\end{array}} \right](x ,x_3 ,t) = {\mathop{\rm Re}\nolimits} \left( {\left[ {\begin{array}{*{20}c}
   {\rm E}  \\
   {\rm H}  \\
\end{array}} \right](x)e^{i(\beta x_3  - \omega t)}} \right)
\end{array}
\end{equation}
are called the eigenwaves of the waveguide. Here, positive $\omega$ is the radian
frequency, $\beta$  is the propagation constant, $\rm E$ and $\rm H$ are complex
amplitudes of $\cal E$ and $\cal H$, $x=(x_1,x_2)$.

In forward eigenvalue problems the permittivity is known and it is necessary to calculate
longitudinal wavenumbers $k=\omega\sqrt{\varepsilon_0\mu_0}$ and propagation constants
$\beta$ such that there exist eigenmodes. The eigenmodes have to satisfy transparency
conditions at the boundaries  and satisfy an condition at infinity. In inverse problems
considering in this work it is necessary to reconstruct the unknown
permit\-tivi\-ty~$\varepsilon$ by some information on natural eigenwaves which exist for
some eigenvalues~$k$ and~$\beta$.

The domain  $\Omega_e=\{x\in\mathbb{R}^2:r>r_n\}$  is unbounded. Therefore, it is
necessary to formulate a condition at infinity for  complex amplitudes~$\mathrm{E}$
and~$\mathrm{H}$ of eigenmodes.  Let us confine ourselves to the investigation of the
surface modes only. The propagation constants $\beta$ of  surface modes are real and
belong to the interval~$G=(k\varepsilon_e,k\varepsilon_+)$. The amplitudes of  surface
modes satisfy to the following condition:
\begin{equation}\label{kem_surf_waves_inf_cond}
\left[ {\begin{array}{*{20}c} {\rm E} \\ {\rm H} \\
\end{array}}\right] = e ^{ {-\sigma r}} O\left(
{\frac{1}{\sqrt{r}}}\right) ,\quad r\rightarrow \infty.
\end{equation}
Here,  $\sigma=\sqrt{\beta ^{2}-k^2\varepsilon_e}>0$ is the transverse wavenumber in the
cladding.

Denote by~$\chi_l=\sqrt{k^2\varepsilon_l-\beta^2}$ the transverse wavenumbers in the
waveguide's layers. Under the weakly guidance approximation~\cite{Snider_Love_1983} the
original problem is re\-du\-ced to the calculation of numbers~$k$ and~$\beta$ such that
there exist nontrivial solutions $u=\rm H_1=H_2$ of Helmholtz equations
\begin{equation}\label{1}
\Delta u + \chi_l ^2  u = 0,\quad  x \in \Omega_l,\quad l=1,\,2,\,\ldots,\, n,
\end{equation}
\begin{equation}\label{2}
\Delta u - \sigma^2  u = 0,\quad  x \in \Omega_e,
\end{equation}
which satisfy  the transparency conditions
\begin{equation}\label{3}
u^{+} = u^{-}, \quad \frac{\partial u^{+}}{\partial r} = \frac{\partial u^{-}}{\partial
r}, \quad x\in\gamma_l,\quad l=1,\,2,\,\ldots,\, n.
\end{equation}
Here $\Omega_1=\{x\in\mathbb{R}^2:r<r_1\}$,
$\Omega_l=\{x\in\mathbb{R}^2:r_{l-1}<r<r_l\}$, $\gamma_l=\{x\in\mathbb{R}^2:r=r_l\}$,
$l=1,\,2,\,\ldots,\, n$;  $u^{-}$ (respectively, $u^{+}$) is the limit value of a
function~$u$ from the interior (respectively, the exterior) of $\gamma_l$.

Let us calculate nontrivial solutions $u$ of problem~(\ref{1})--(\ref{3}) in the space of
continuous and continuously differentiable in $\overline{\Omega}_e$ and
$\overline{\Omega}_l$, $l=1,\,2,\,\ldots,\, n$, and twice continuously differentiable in
$\Omega_e$ and $\Omega_l$, $l=1,\,2,\,\ldots,\, n$, functions, satisfying  condition
\begin{equation}\label{4}
u = e ^{ {-\sigma r}} O\left( {\frac{1}{\sqrt{r}}}\right) ,\quad r\rightarrow \infty.
\end{equation}
Denote by $U$ described functional space. We solve problem (\ref{1})--(\ref{3}) by the
method of separation of variables using polar coordinates $(r,\varphi)$ and look for the
function $u\in U\setminus \{0\}$ in the form
\begin{equation}\label{u_repr_RP}
    u(r,\varphi)=R(r)\Phi(\varphi).
\end{equation}
Then
$$
\Phi(\varphi)= \left\{\begin{array}{c}
                \cos m\varphi \\
                \sin m\varphi \\
              \end{array}\right.,\quad m=0,\,1,\,2,\,\ldots,
$$
$$
R(r)=\left\{\begin{array}{c}
       aJ_m(\chi_1r),\quad r<r_1 \\
       b_lJ_m(\chi_lr)+c_lH_m^{(1)}(\chi_lr),\quad r_{l-1}<r<r_{l},\quad l=2,\,\ldots,\, n        \\
       dK_m(\sigma r), \quad r> r_n. \\
     \end{array}\right.,
$$
here $J_m(z)$ is the Bessel function, $H_m^{(1)}(z)$ is the Hankel function of the first
kind, $K_m(z)$ is the Macdonald function. Unknown coefficients $a$, $b_l$, $c_l$, and~$d$
satisfy the following homogeneous system of linear algebraic equations:
\begin{equation}\label{main_SLAE}
    \left(%
\begin{array}{cccc}
  A_{11}(k,\beta,\varepsilon_1,\varepsilon_2) & 0                           & \ldots & 0 \\
  0                           & A_{22}(k,\beta,\varepsilon_2,\varepsilon_3) & \ldots & 0 \\
  \ldots                      & \ldots                      & \ldots &  \ldots \\
  0                           & 0                           & \ldots & A_{nn}(k,\beta,\varepsilon_n,\varepsilon_e) \\
\end{array}%
\right)
\left(%
\begin{array}{c}
  X_1 \\
  X_2 \\
  \ldots \\
  X_n \\
\end{array}%
\right)=0.
\end{equation}
Here
$$
A_{11}(k,\beta,\varepsilon_1,\varepsilon_2)=\left(%
\begin{array}{ccc}
  J_m(\chi_1 r_1) & -J_m(\chi_2 r_1) & -H_m^{(1)}(\chi_2 r_1) \\
  \chi_1 J'_m(\chi_1 r_1) & -\chi_2J'_m(\chi_2 r_1) & -\chi_2 H_m^{(1)'}(\chi_2 r_1) \\
\end{array}%
\right),
$$
$$
A_{nn}(k,\beta,\varepsilon_n,\varepsilon_e)=\left(%
\begin{array}{ccc}
  J_m(\chi_n r_n)         & H_m^{(1)}(\chi_n r_n)         & -K_m(\sigma r_n) \\
  \chi_n J'_m(\chi_n r_n) & \chi_n H_m^{(1)'}(\chi_n r_n) & -\sigma K'_m(\sigma r_n) \\
\end{array}%
\right),
$$
$$
A_{ll}(k,\beta,\varepsilon_l,\varepsilon_{l+1})=
$$
$$
=\left(%
\begin{array}{cccc}
  J_m(\chi_l r_l)         & H_m^{(1)}(\chi_l r_l)         &-J_m(\chi_{l+1} r_l)            & -H_m^{(1)}(\chi_{l+1} r_l) \\
  \chi_l J'_m(\chi_l r_l) & \chi_l H_m^{(1)'}(\chi_l r_l) &-\chi_{l+1}J'_m(\chi_{l+1} r_l) & -\chi_{l+1} H_m^{(1)'}(\chi_{l+1} r_l) \\
\end{array}%
\right),
$$
$$
X_1=\left(%
\begin{array}{c}
  a \\
  b_2 \\
  c_2 \\
\end{array}%
\right),\ X_n=\left(%
\begin{array}{c}
  b_n \\
  c_n \\
  d \\
\end{array}%
\right),\ X_l=\left(%
\begin{array}{c}
  b_l \\
  c_l \\
  b_{l+1} \\
  c_{l+1} \\
\end{array}%
\right),
$$
where $l=2,\,3,\,\ldots,\, n-1$. System \eqref{main_SLAE} has a nontrivial solution if
and only if the determinant of its matrix is equal to zero:
\begin{equation}\label{main_char_eq}
\det (A(k,\beta,\varepsilon)) =0.
\end{equation}
The last condition in the theory of optical waveguides is called the characteristic
equation.

\section{Ill-posed problems}

In this section on the base of characteristic equation \eqref{main_char_eq} we formulate
 problems of reconstruction of the unknown permit\-tivi\-ty~$\varepsilon$ by some
information on the fundamental eigenwaves which exist for some~$k$ and~$\beta$. In the
case of the fundamental eigenwave the permittivity $\varepsilon$ are connected with each
value of the propagation constant $\beta$ and each value of the longitudinal
wavenumber~$k$ by
 characteristic equation \eqref{main_char_eq} with $m=0$:
\begin{equation}\label{kem_char_eq}
   \det (A(k,\beta,\varepsilon)) = \left|%
\begin{array}{cccc}
  A_{11}(k,\beta,\varepsilon_1,\varepsilon_2) & 0                           & \ldots & 0 \\
  0                           & A_{22}(k,\beta,\varepsilon_2,\varepsilon_3) & \ldots & 0 \\
  \ldots                      & \ldots                      & \ldots &  \ldots \\
  0                           & 0                           & \ldots & A_{nn}(k,\beta,\varepsilon_n,\varepsilon_e) \\
\end{array}%
\right| = 0.
\end{equation}
Here
$$
A_{11}(k,\beta,\varepsilon_1,\varepsilon_2)=\left(%
\begin{array}{ccc}
  J_0(\chi_1 r_1) & -J_0(\chi_2 r_1) & -H_0^{(1)}(\chi_2 r_1) \\
  \chi_1 J'_0(\chi_1 r_1) & -\chi_2J'_0(\chi_2 r_1) & -\chi_2 H_0^{(1)'}(\chi_2 r_1) \\
\end{array}%
\right),
$$
$$
A_{nn}(k,\beta,\varepsilon_n,\varepsilon_e)=\left(%
\begin{array}{ccc}
  J_0(\chi_n r_n)         & H_0^{(1)}(\chi_n r_n)         & -K_0(\sigma r_n) \\
  \chi_n J'_0(\chi_n r_n) & \chi_n H_0^{(1)'}(\chi_n r_n) & -\sigma K'_0(\sigma r_n) \\
\end{array}%
\right),
$$
$$
A_{ll}(k,\beta,\varepsilon_l,\varepsilon_{l+1})=
$$
$$
=\left(%
\begin{array}{cccc}
  J_0(\chi_l r_l)         & H_0^{(1)}(\chi_l r_l)         &-J_0(\chi_{l+1} r_l)            & -H_0^{(1)}(\chi_{l+1} r_l) \\
  \chi_l J'_0(\chi_l r_l) & \chi_l H_0^{(1)'}(\chi_l r_l) &-\chi_{l+1}J'_0(\chi_{l+1} r_l) & -\chi_{l+1} H_0^{(1)'}(\chi_{l+1} r_l) \\
\end{array}%
\right),
$$
where $l=2,\,3,\,\ldots,\, n-1$.

We suppose that the permittivities $\varepsilon_1$, $\varepsilon_2$, ... $\varepsilon_n$,
$\varepsilon_e$ are real and satisfy  conditions \eqref{kem_perm_cond_1}. Denote by
$\varepsilon=(\varepsilon_1,\varepsilon_2,\ldots,\varepsilon_n,\varepsilon_e)\in
\mathbb{R}^{n+1}$ the vector of the permittivities of the waveguide. Then we can write
condition \eqref{kem_perm_cond_1} in the matrix form:
\begin{equation}\label{kem_perm_cond_2}
C\varepsilon \leq q,
\end{equation}
$$
C=
\left(%
\begin{array}{ccccc}
  \phantom{-}0      & 0      & \ldots            & \phantom{-}0      & -1 \\
  -1                & 0      & \ldots            & \phantom{-}0      & \phantom{-}1 \\
  \phantom{-}\ldots & \ldots & \phantom{-}\ldots & \ldots            & \phantom{-}\ldots \\
  \phantom{-}0      & 0      & \ldots            & -1                & \phantom{-}1 \\
\end{array}%
\right),\quad q=
\left(%
\begin{array}{c}
  -1 \\
  -\mu \\
  -\mu \\
  -\mu \\
\end{array}%
\right),
$$
where $\mu$ is a small positive number.

Suppose that we know $m\geq n+1$ pairs of values of  the longitudinal wavenumber and the propagation
constant of the fundamental eigenwave of the waveguide:~$(k_i,\beta_i)$, $i=1,\,2,\,\ldots,\,m$. Let us
introduce $m$ functions $f_i(\varepsilon)$ of the variable $\varepsilon$ by the following way:
%$$
%f_i(\varepsilon)=|\det(A(k_i,\beta_i,\varepsilon))|, \quad i=1,\,2,\,\ldots,\,m,
%$$
%or
$$
f_i(\varepsilon)=1/\textrm{cond}(A(k_i,\beta_i,\varepsilon)), \quad i=1,\,2,\,\ldots,\,m.
$$
Here $A$ is the matrix defined in \eqref{kem_char_eq}. By $\textrm{cond}$ we denote the condition number
of $A$. Clearly, all $f_i$ are non-negative and if $\varepsilon$ satisfies characteristic
equation~\eqref{kem_char_eq}, then $f_i(\varepsilon)=0$, $i=1,\,2,\,\ldots,\,m$. Let us introduce the
nonlinear operator $F:\mathbb{R}^{n+1}\to \mathbb{R}^m$ by the formula:
$$
F(\varepsilon)=(f_1(\varepsilon),f_2(\varepsilon),\ldots,f_m(\varepsilon)).
$$

It seems natural to find the vector $\varepsilon$ of the permittivities as a solution of the problem
\begin{equation}\label{kem_ill_prob_1}
    F(\varepsilon)=0,\quad C\varepsilon\leq q,
\end{equation}
it also seems useful  to study along with \eqref{kem_ill_prob_1} the minimization problem
\begin{equation}\label{kem_ill_prob_2}
    \min\limits_{C\varepsilon \leq q}\frac{1}{2} ||F(\varepsilon)||^2_{\mathbb{R}^m},
\end{equation}
but this way is not correct. Problems \eqref{kem_ill_prob_1} and \eqref{kem_ill_prob_2} are not
equivalent; each solution of \eqref{kem_ill_prob_1} is evidently a global minimizer of problem
\eqref{kem_ill_prob_2} while a solution of  \eqref{kem_ill_prob_2} do not necessarily satisfy
\eqref{kem_ill_prob_1}. Moreover, problems \eqref{kem_ill_prob_1} and \eqref{kem_ill_prob_2} are
ill-posed (see, for example \cite{BK}). The ill-posedness of these problems means that analyzing
problems close in some sense to \eqref{kem_ill_prob_1} and \eqref{kem_ill_prob_2}, we can not guarantee
that solutions to these perturbed problems are close to corresponding solutions of the original ones.

\section{The Tikhonov functional}

In real-life applications for each fixed longitudinal wavenumber the propagation constant
of the fundamental eigenwave of the waveguide is measured by physical experiments. Denote
by $({k}_i,\widetilde{\beta}_i)$, $i=1,\,2,\,\ldots,\,m$, these values. Denote by
$\widetilde{F}:\mathbb{R}^{n+1}\to \mathbb{R}^m$ the perturbed operator, where
$$
\widetilde{F}(\varepsilon)=(\widetilde{f}_1(\varepsilon),\widetilde{f}_2(\varepsilon),\ldots,\widetilde{f}_m(\varepsilon)),
$$
$$
\widetilde{f}_i(\varepsilon)=1/\textrm{cond}(A({k}_i,\widetilde{\beta}_i,\varepsilon)),
\quad i=1,\,2,\,\ldots,\,m.
$$
Therefore instead of original problem \eqref{kem_ill_prob_1} we have perturbed problem
\begin{equation}\label{kem_1}
    \widetilde{F}(\varepsilon)=0,\quad C\varepsilon\leq q.
\end{equation}
We suppose that the perturbation is small, namely,
$$
||\widetilde{F}(\varepsilon_*)||_{\mathbb{R}^m}\leq\delta.
$$
Here the small parameter $\delta \in (0,1)$ characterizes the level of the error in data,
$\varepsilon_*$ is the exact solution of non-perturbed problem \eqref{kem_ill_prob_1}.

As we have seen,  problem \eqref{kem_ill_prob_1} is a classical ill-posed problem \cite{TLY}. Thus, we
assume that there exists the exact solution~$\varepsilon_*$ to our problem \eqref{kem_ill_prob_1} but we
never will get this solution in computations. Because of that we call by the \emph{regularized solution}
$\varepsilon_{\alpha}$ some approximation of the unknown exact solution~$\varepsilon_*$ which is
satisfied to the requirements of closeness to the exact solution $\varepsilon_*$ and stability with
respect to the small errors $\delta$.

We use Tikhonov regularization algorithm (see \cite{TLY}) which is based on the minimization of the
Tikhonov functional. Thus, to find regularized solution~$\varepsilon_{\alpha}$ of problem (\ref{kem_1}),
we minimize the Tikhonov regularization functional
\begin{equation}
M_{\alpha}(\varepsilon) =\frac{1}{2}\left\Vert \widetilde{F}(\varepsilon) \right\Vert
_{\mathbb{R}^m}^{2}+\frac{\alpha}{2}\left\Vert \varepsilon-
\varepsilon_{0}\right\Vert_{\mathbb{R}^{n+1}}^{2},\quad C\varepsilon\leq q,
 \label{tikh5}
\end{equation}
\begin{equation*}
M_{\alpha}: \mathbb{R}^{n+1}\rightarrow \mathbb{R},\quad \varepsilon_{0}\in \mathbb{R}^{n+1},
\end{equation*}
where $\alpha =\alpha \left( \delta \right) >0$ is a small regularization parameter. The choice of  the
point $\varepsilon_{0}$ and the regularization parameter $\alpha$ depends on the concrete minimization
problem. This question will be investigated  later by numerical experiments. Usually $\varepsilon_{0}$
is a good first approximation for the exact solution~$\varepsilon_*.$

It follows from \cite{TLY} that an algorithm for solution of the equation (\ref{kem_1})
which is based on the minimization of the Tikhonov functional (\ref{tikh5}) is the
regularization algorithm, and the element $\varepsilon_{\alpha} \in \mathbb{R}^{n+1}$
where the functional (\ref{tikh5}) reaches its minimum is the regularized solution.

In our theoretical investigations below we need reformulate results of
\cite{BKK, BOOK} for the case of our \textbf{IP}. In this section
below $||\cdot||$ denotes $\mathbb{R}^{n+1}$ norm.  Let $H_1$ be the
finite dimensional linear space.  Let $Y$ be  the set of admissible
parameters for $\varepsilon \in \mathbb{R}^{n+1}$ defined in
(\ref{kem_perm_cond_2}) and let us define by $Y_1 := Y \cap H_1$ with
$G := \bar{Y}_1$.  Now the operator $\widetilde{F}:G \to H_2$
corresponds to the operator in the Tikhonov functional (\ref{tikh5}).

We now assume that the operator $\widetilde{F}(\varepsilon)$ defined in (\ref{kem_1}) is one-to-one
and denote by   $V_{r}\left(\varepsilon\right)$  neighborhood of the radius $r$ of
$\varepsilon$  such that
\begin{equation}\label{neigh}
V_{r}\left(\varepsilon\right) =\left\{\varepsilon^{\prime }\in H_{1}:\left\|
\varepsilon^{\prime } - \varepsilon \right\|   < r,~~ \forall r >0 ~~\forall \varepsilon
\in H_{1}\right\} .
\end{equation}
We also  make common assumptions, see for details \cite{BKS, tikhonov}, that the operator
$\widetilde{F}$ has the Lipschitz continuous Frech\'{e}t derivative $\widetilde{F}^{\prime
}(\varepsilon) $ for $\varepsilon \in V_{1}(\varepsilon^{\ast}),$  such that there exist
constants  $N_{1},N_{2} > 0$
\begin{equation}
\left\| \widetilde{F}^{\prime }(\varepsilon) \right\| \leq N_{1},\left\| \widetilde{F}^{\prime
}(\varepsilon_1) - \widetilde{F}^{\prime }(\varepsilon_2) \right\| \leq N_{2}\left\|
\varepsilon_1 - \varepsilon_2 \right\| ,\forall \varepsilon_1, \varepsilon_2 \in
V_{1}\left(\varepsilon^{\ast }\right) .  \label{2.7}
\end{equation}

Similarly with \cite{BKK} we choose the constant $D= D\left( N_{1},N_{2}\right)
=const.>0$ such that
\begin{equation}
| M_{\alpha}^{\prime }(\varepsilon_1) -  M_{\alpha}^{\prime }(\varepsilon_2) | \leq
D\left\| \varepsilon_1 - \varepsilon_2\right\| ,\forall \varepsilon_1,\varepsilon_2 \in
V_{1}(\varepsilon^*). \label{2.10}
\end{equation}
Through the paper as in \cite{BKK}  we assume that
\begin{eqnarray}
\left\| \varepsilon_0 - \varepsilon^{\ast }\right\| &\leq &\delta ^{\xi},~\xi =const.\in \left( 0,1\right) ,  \label{2.11} \\
\alpha &=&\delta ^{\zeta}, \zeta=const.\in ( 0, \min(\xi, 2 - 2\xi)),  \label{2.12}
\end{eqnarray}
where  $\alpha$ is the regularization parameter. Equation (\ref{2.11}) means that we
assume that all initial guess $\varepsilon_0$  in the Tikhonov functional is located in a
sufficiently small neighborhood $V_{\delta ^{\xi}}(\varepsilon^*)$ of the exact solution
$\varepsilon^*$. From Lemmata 2.1 and 3.2 of \cite{BKK} follows that conditions
(\ref{2.12}) ensures that  $(\varepsilon^{\ast }, \varepsilon_0)$ belong to an
appropriate neighborhood of the regularized solution of the Tikhonov functiona.

Below we reformulate Theorem 1.9.1.2 of \cite{BOOK} for the case of our Tikhonov
functional. Different proofs of this theorem can be found in \cite{BOOK} and in
\cite{BKK} and are straightly applied to our   case.

\textbf{Theorem 1} \emph{Let }$\Omega \subset \mathbb{R}^{3}$ \emph{\ be a convex bounded
domain with the boundary }$\partial \Omega \in C^{3}.$   \emph{ Assume that there exists
the exact
    solution }$\varepsilon^* \in G$\emph{\ of the equation }$\widetilde{F}(\varepsilon^*)
  =0$ \emph{\ for the case of the exact data }$ (k_i^\ast, \beta_i^\ast)$.\emph{  Let regularization parameter } $\alpha$
  \emph{ in (\ref{tikh5}) is such that }
\begin{equation*}
\alpha = \alpha(\delta) =\delta ^{2\nu},~~\nu  =const.\in \left( 0,\frac{1}{4}\right),
~\quad \forall \delta \in \left( 0,1\right).
\end{equation*}
\emph{Let } $\varepsilon_0$ \emph{ satisfies  conditions
  (\ref{2.11}). Then the Tikhonov functional  (17)
    is strongly convex in the neighborhood }$V_{\alpha} \left( \delta
    \right) (\varepsilon^*) $ \emph{\ with the strong convexity constant }$\gamma =  \alpha.$
 \emph{The strong convexity property can be also written as}
\begin{equation}
\left\Vert \varepsilon_{1} - \varepsilon_{2}\right\Vert ^{2}\leq \frac{2}{\delta ^{2 \nu
}}\left( M_{\alpha}'(\varepsilon_1) - M_{\alpha}'(\varepsilon_2), \varepsilon_1 -
\varepsilon_{2}\right), \text{ }\forall \varepsilon_1, \varepsilon_2 \in H_{1},
\label{4.249}
\end{equation}
where $(\cdot, \cdot)$ is a scalar product. \emph{Next, there exists
  the unique regularized solution } $\varepsilon_{\alpha}$ \emph{ of the
  functional (17) and this solution } $\varepsilon_{\alpha} \in
V_{\delta ^{3\nu }/3}(\varepsilon^*).$\emph{\ The
  gradient method of the minimization of the functional
  (\ref{tikh5}) which starts at }$\varepsilon_0$\emph{\ converges to the regularized solution of this
  functional. Furthermore,}
\begin{equation}\label{accur}
\left\Vert \varepsilon_{\alpha}  - \varepsilon^* \right\Vert \leq \theta \left\Vert
\varepsilon_0 - \varepsilon^* \right\Vert, ~~\theta \in (0,1).
\end{equation}

The property(\ref{accur}) means that the regularized solution of the Tikhonov functional
(\ref{tikh5}) provides a better accuracy than the initial guess $\varepsilon_0$  if it satisfies
condition (\ref{2.11}).

\section{Numerical experiments}

We minimized the Tikhonov functional (\ref{tikh5}) using the GlobalSearch Algorithm of
the GlobalSearch object in MATLAB in the case of the one-layered waveguide. The exact
values of parameters are chosen as follows:
 $\varepsilon_i=\varepsilon_1=2.383936$ (quartz) and
$\varepsilon_e=2.21235876$ (optical glass).
 The exact solutions $(k_i^2,\beta_i^2)$ of
the forward spectral problem for this waveguide are well known (see, for example,
\cite{Spiridonov_Karchevskiy_DD_2013}). We constructed the Tikhonov functional using five
pairs of eigenvalues corresponding to the fundamental eigenwave, see Table 1.

\begin{table}[tbp]
{\footnotesize Table 1. \emph{Five
pairs of eigenvalues corresponding to the fundamental eigenwave  which are used in numerical tests.}}  \par
\vspace{2mm}
\begin{center}
\begin{tabular}{|c|c|c|}
  \hline
  $i$ & $k_i^2$ & $\beta_i^2$ \\
  \hline
  1 & 02.4 & 05.30978783787819 \\
  2 & 04.8 & 10.63618822212100 \\
  3 & 07.2 & 16.02129849868130 \\
  4 & 09.6 & 21.46863534179760 \\
  5 & 12.0& 26.96464966481510 \\
  \hline
\end{tabular}
\end{center}
\end{table}

It follows from the results of \cite{Spiridonov_Karchevskiy_DD_2013} that for the wide
range of frequencies the following formula gives a very good approximation from above to
the permittivity~$\varepsilon_e$ of the cladding of the waveguide:
$$
\varepsilon_e\approx {\beta_1^2}/{k_1^2}.
$$
Using this formula and condition \eqref{kem_perm_cond_2}, we took the first approximation
to the vector of permittivities $\varepsilon=(\varepsilon_i,\,\varepsilon_e)$ as
$$
\varepsilon_0=(\varepsilon_{0,i},\varepsilon_{0,e})=({\beta_1^2}/{k_1^2},\,
{\beta_1^2}/{k_1^2})=(2.21241159911591,\,2.21241159911591).
$$

In our computations by analogy with \cite{Karchecvskii_Spiridonov_Beilina_SPMS_2015} we
introduced a randomly distributed noise in the propagation constants as
$$
\widetilde{\beta_i}=\beta_i(1+p\alpha),\quad i=1,\,2,\,\ldots\,5,
$$
where $\beta_i$ are the exact measured propagation constants, $\alpha\in(-1,1)$ are
randomly distributed numbers, and~$p$ is the noise level. In our computations we used
$p=0.05$ and thus, the noise level was~$5\%$. In the numerical experiments we accounted
this noise in perturbed operator \eqref{kem_1}  and also in the perturbed initial
approximations
$$
\widetilde{\varepsilon}_0=(\widetilde{\varepsilon}_{0,i},
\widetilde{\varepsilon}_{0,e})=({\widetilde{\beta}_1^2}/{k_1^2},\,
{\widetilde{\beta}_1^2}/{k_1^2})
$$
that we used in \eqref{tikh5} instead of $\varepsilon_0$. The regularization parameter
was $\alpha=0.01$.

\begin{figure}[b!]\centering
  \includegraphics[width=1\textwidth]{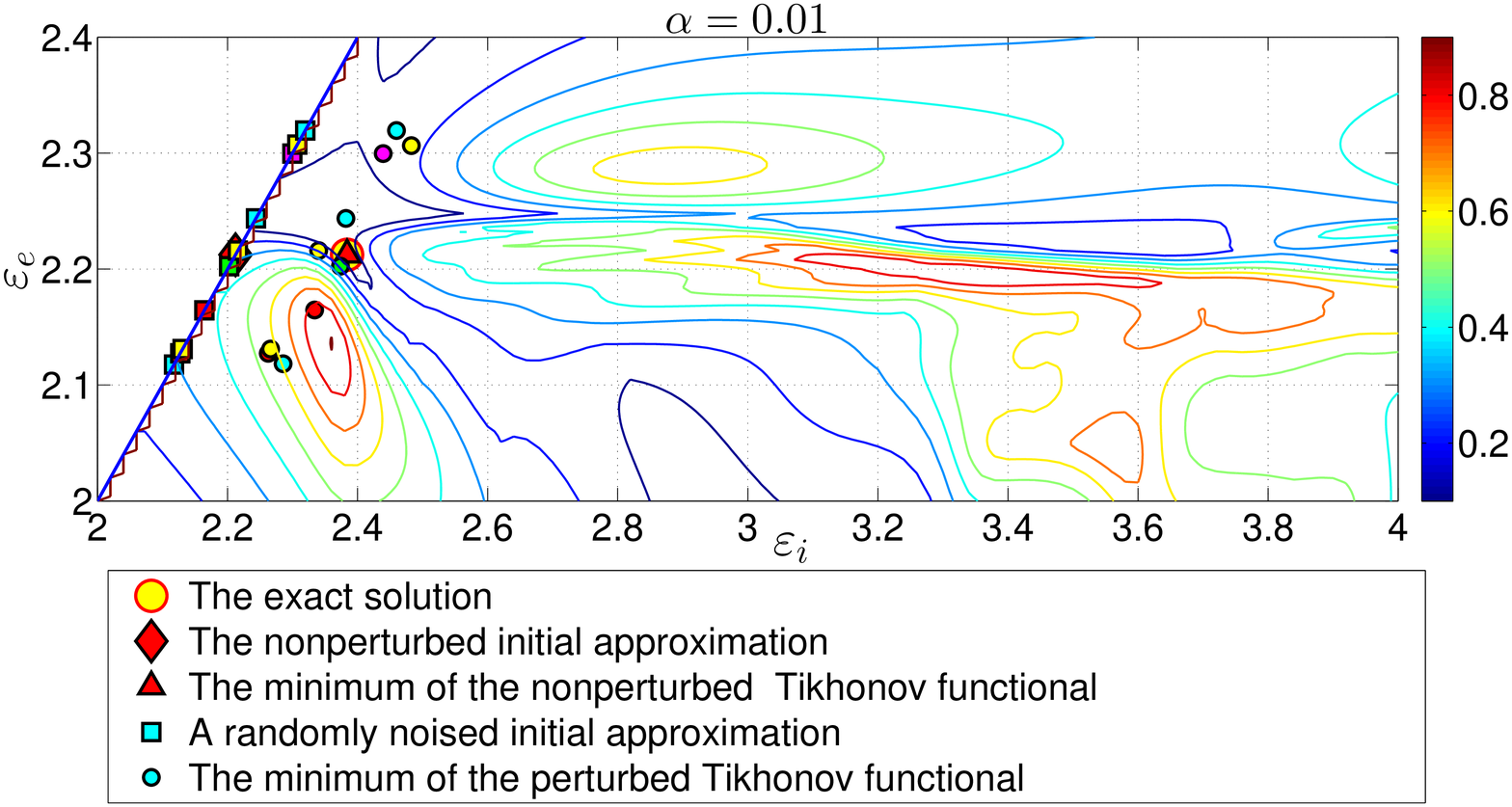}
   \caption{{\footnotesize \emph{The results of numerical minimization of the nonperturbed and perturbed
  (by the randomly noised propagation constants  $\widetilde{\beta_i}$ with the $5\%$ noise level) Tikhonov
  functional.  The background is the pattern of the nonperturbed
  functional.}}}\label{plot_2}
\end{figure}

\begin{figure}[h!]\centering
  \includegraphics[width=1\textwidth]{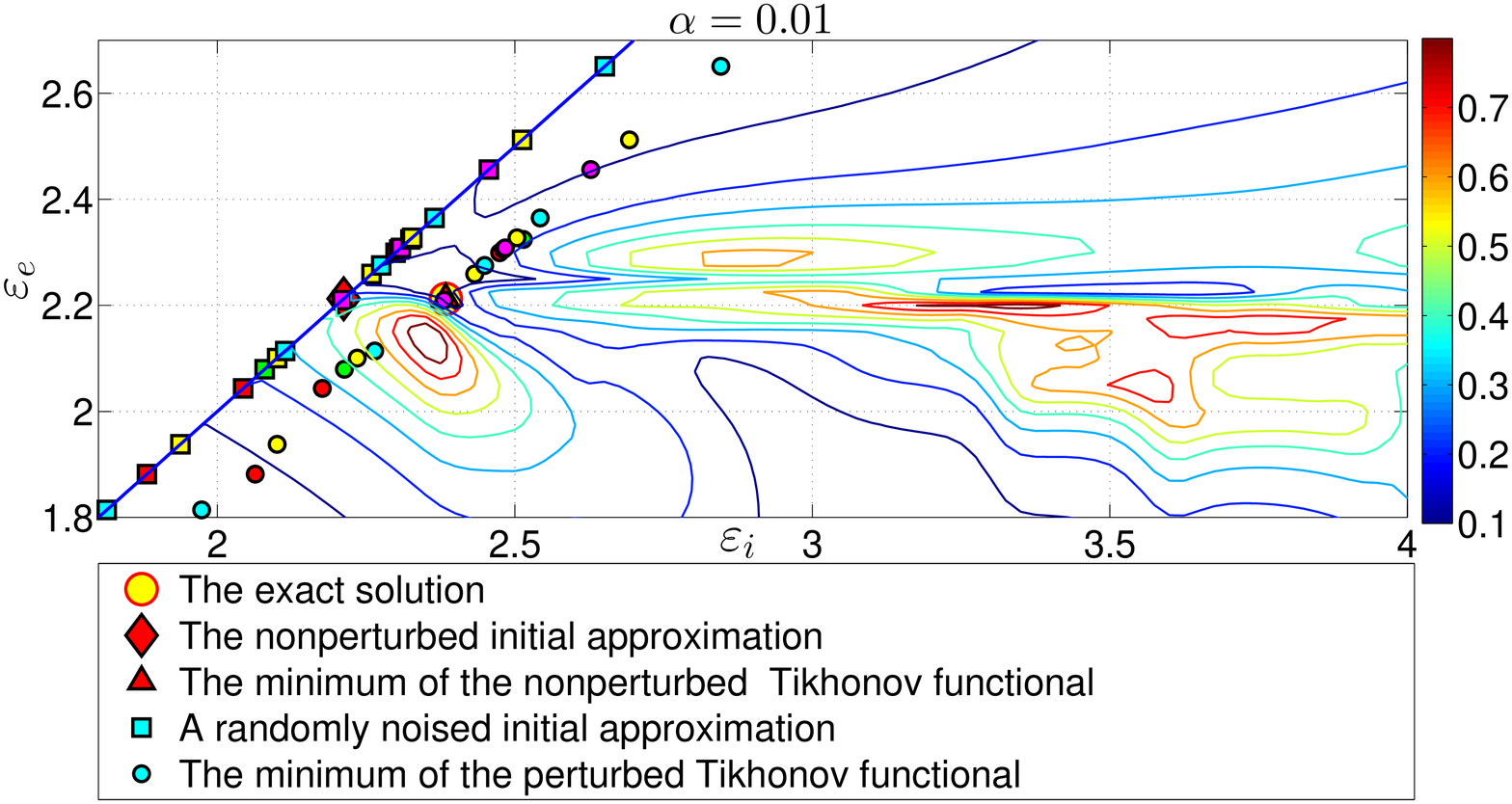}
  \caption{ {\footnotesize \emph{ The results of numerical minimization of the nonperturbed and perturbed
  (by the randomly noised propagation constants  $\widetilde{\beta_i}$ with the $20\%$ noise level) Tikhonov
  functional.  The background is the pattern of the nonperturbed
  functional.}}}\label{plot_3}
\end{figure}

Some numerical results of reconstruction of the vector~$\varepsilon$ of the permittivies
are presented at Figure \ref{plot_2}. The exact value $\varepsilon$  is marked  at
Figure \,\ref{plot_2} by the big yellow circle. The approximated value
$$\varepsilon_{\alpha}=(2.38393603911088,\, 2.21235872913944)$$ of $\varepsilon$ for the
noise-free data is marked by the red triangle. The background of this figure is the
pattern of the nonperturbed Tikhonov functional. Approximated values
$\widetilde{\varepsilon}_{\alpha}=(\widetilde{\varepsilon}_{\alpha, i},
\widetilde{\varepsilon}_{\alpha, e})$ of $\varepsilon$ for randomly distributed noise
$\widetilde{\beta_i}$ with the $5\%$ noise level are marked by the colored circles. The
nonperturbed initial approximation~${\varepsilon}_0$ is marked by the red rhomb, and the
perturbed initial approximations~$\widetilde{\varepsilon}_0$ are marked by the colored
squares. The Table 2 presents results of the reconstruction for  different initial guesses   with random noise level $\sigma=5\%$ in data.
 Using the Figure \ref{plot_3} and  the Table 2 we observe that the
approximate solutions were stable even for the randomly noised~$\widetilde{\beta_i}$. Using the  Table 2 we also observe that the computed
relative error $e = {\Vert\varepsilon
-\widetilde{\varepsilon}_{\alpha}\Vert_{\mathbb{R}^{2}}
}/{\Vert\varepsilon\Vert_{\mathbb{R}^{2}} }$ is on the interval $[0,\,0.05]$, and
 thus, approximated values $\widetilde{\varepsilon}_{\alpha}$ differs from the exact values of
$\varepsilon$ not more than $5\%$.

\medskip
\begin{table}[tbp]
{\footnotesize Table 2. \emph{Computational results of the
    reconstructions $\widetilde{\varepsilon}_{\alpha,i}, \widetilde{\varepsilon}_{\alpha,e}$ together with computational errors $e$  for different initial guesses  $\widetilde{\varepsilon}_{0,i} = \widetilde{\varepsilon}_{0,e}$. Noise in data
    is $\sigma=5\%$.}}  \par
\vspace{2mm}
\begin{tabular}{|c|c|c|c|}
  \hline
  $\widetilde{\varepsilon}_{0,i} = \widetilde{\varepsilon}_{0,e}$&
  $\widetilde{\varepsilon}_{\alpha,i}$ & $\widetilde{\varepsilon}_{\alpha,e}$ & $e$ \\
  \hline
  2.11768385759666 & 2.28552927271664 & 2.11830350061034 & 0.04186 \\
  2.12734296884777 & 2.26254671102130 & 2.15032840068311 & 0.04192 \\
  2.13098610290990 & 2.26628839765012 & 2.15401089492431 & 0.04038 \\
  2.16434434702641 & 2.33397022113028 & 2.16467608183591 & 0.02124 \\
  2.20248954651734 & 2.37362289984480 & 2.20251745727489 & 0.00439 \\
  2.21597381604893 & 2.34013490935047 & 2.25067212773741 & 0.01790 \\
  2.24374276867624 & 2.38204127941254 & 2.26798587326534 & 0.01712 \\
  2.29968977385246 & 2.43945637236437 & 2.32453737220619 & 0.03849 \\
  2.30737609305936 & 2.48261925050832 & 2.30651724858946 & 0.04194 \\
  2.31946862194614 & 2.45975978096066 & 2.34452993693523 & 0.04686 \\
  \hline
\end{tabular}
\end{table}

\medskip

Analogous results we obtained for the $20\%$ noise level. We present
reconstruction results in Figure \ref{plot_3} and in the Table 3. Using the Table 3 we see,  that  the relative error $e$ is located  on the interval $[0,\,0.20]$, and
approximated values $\widetilde{\varepsilon}_{\alpha}$ differs from
the exact values of $\varepsilon$ not more than $20\%$.

\medskip
\begin{table}[tbp]
{\footnotesize Table 3. \emph{Computational results of the
    reconstructions $\widetilde{\varepsilon}_{\alpha,i}, \widetilde{\varepsilon}_{\alpha,e}$ together with computational errors $e$  for different initial guesses  $\widetilde{\varepsilon}_{0,i} = \widetilde{\varepsilon}_{0,e}$. Noise in data
    is $\sigma=20\%$.}}  \par
\begin{tabular}{|c|c|c|c|}
  \hline
  $\widetilde{\varepsilon}_{0,i} = \widetilde{\varepsilon}_{0,e}$&
  $\widetilde{\varepsilon}_{\alpha,i}$ & $\widetilde{\varepsilon}_{\alpha,e}$ & $e$ \\
  \hline
   1.81417751127505 &  1.97381469634454 &  1.81480974434505 & 0.17563  \\
   1.88178880974403 &  2.06363445374894 &  1.87910974743464 & 0.14212  \\
   1.93807256082554 &  2.10049008805235 &  1.93817330919554 & 0.12126  \\
   2.04345147329207 &  2.17638931794059 &  2.06553048164261 & 0.07817  \\
   2.07966690316895 &  2.21358678746843 &  2.10213721049311 & 0.06239  \\
   2.10053329208657 &  2.23501858364500 &  2.12322904231478 & 0.05337  \\
   2.11410009015371 &  2.26444198310842 &  2.14720323360852 & 0.04185  \\
   2.20972436445250 &  2.38114288964767 &  2.20969348288490 & 0.00119  \\
   2.25920990213314 &  2.43257430919639 &  2.25876527535809 & 0.02067  \\
   2.27540878243121 &  2.44940844145173 &  2.27482428589175 & 0.02783  \\
   2.29950015590045 &  2.47444275383217 &  2.29870392359414 & 0.03847  \\
   2.30564544562698 &  2.48082751592143 &  2.30479520038292 & 0.04118  \\
   2.30861157840774 &  2.48391025002180 &  2.30773406680037 & 0.04249  \\
   2.32443445081245 &  2.51439918461942 &  2.36083107170200 & 0.06078  \\
   2.32757009410069 &  2.50359310544044 &  2.32653791258360 & 0.05086  \\
   2.36509793984544 &  2.54259737446660 &  2.36370206817162 & 0.06742  \\
   2.45636539769611 &  2.62761067485400 &  2.45337536941920 & 0.10539  \\
   2.51223352049175 &  2.69257246813461 &  2.51223352049175 & 0.13232  \\
   2.65120325027458 &  2.84608346134752 &  2.64927524365834 & 0.19555  \\
  \hline
\end{tabular}
\end{table}

\medskip

Using Tables 2, 3 and Figures \ref{plot_2}, \ref{plot_3} we conclude that the
optimization algorithm with the proposed first approximations  for $\varepsilon_0$,
computed using theory of \cite{Frolov_Kartchevskiy_2013, Spiridonov_Karchevskiy_DD_2013},
gives  stable reconstruction of~$\varepsilon$.

\section{Discussion and Conclusion}

In this work, we present new method for the numerical reconstruction
of the variable refractive index of multi-layered circular weakly
guiding dielectric waveguides using the measurements of the
propagation constants of their eigenwaves.  The method is new in the
sense that instead of the conventional measurements of the
time-dependent electrical field we use measurements of the propagation
constants. Such measurements for a multi-layered dielectric waveguide
of arbitrary cross-section can be done in the millimeter range by a
resonance method \cite{Karpenko_1981, Sotsky_1987}.

We  present computational study of the reconstruction of function $\varepsilon$ using
propagation constant measurements. Theorem 1 guarantees convergence of this algorithm in
the case if we have a good first approximation to the function~$\varepsilon$. However,
this is not issue in our case since the analysis of the forward problem developed in
\cite{Frolov_Kartchevskiy_2013, Spiridonov_Karchevskiy_DD_2013} allows us to obtain a
good first approximation to the refractive index of the waveguide, and we use this
initial guess in all our experiments.  We have confirmed our theoretical investigations
by numerical tests, where we have obtained stable reconstruction of the dielectric
permittivity function $\varepsilon$ for random noise level $\sigma=5\%$ and $\sigma=20\%$
in data.

\section*{Acknowledgements}

The research of LB was supported by the Swedish Research Council (VR).

\end{document}